\documentclass[a4paper,11pt]{amsart}
\usepackage[cp1251]{inputenc}
\usepackage[english]{babel}
\usepackage{amsmath,amssymb,euscript,amsthm,amsfonts}
\newtheorem{lem}{Lemma}
\newtheorem{thm}{Theorem}
\theoremstyle{definition}

\theoremstyle{remark}

\theoremstyle{proposition}
\newtheorem{prop}{Proposition}

\begin{document}

\renewcommand{\proofname}{Proof}
\makeatletter \headsep 10 mm \footskip 10 mm


\author{Mathieu Dutour Sikiri\'c}
\address{Mathieu Dutour Sikiri\'c, Rudjer Boskovi\'c Institute, Bijenicka 54, 10000 Zagreb, Croatia, Fax: +385-1-468-0245}
\email{mdsikir@irb.hr}
\thanks{The first author is supported by the Croatian Ministry of Science, Education and Sport under contract 098-0982705-2707}

\author{Viacheslav Grishukhin}
\address{Viacheslav Grishukhin, CEMI Russian Academy of Sciences,
Nakhimovskii prosp.47 117418 Moscow, Russia}
\email{grishuhn@cemi.rssi.ru}

\title{Closure of principal L-type domain and its parallelotopes}

\maketitle

\begin{abstract}
Voronoi defined two polyhedral partitions of the cone of se\-mi\-de\-fi\-nite forms into L-type domains and into perfect domains. Up to equivalence, there is only one domain that is simultaneously perfect and L-type. Voronoi called this domain {\em principal}. We show that closure of the principal domain may be identified with a cone of cut submodular set functions. Parallelotopes of the closed principal domain are zonotopes that are base polyhedra related to graphic unimodular sets of vectors.
\end{abstract}

\section{Introduction}

A parallelotope is a polytope whose translations under a lattice tile the space. A parallelotope is necessarily centrally symmetric with centrally symmetric facets. A zonotope is Minkowski sum of segments. It is known, (see, for example, Proposition 2.2.10 in \cite{BZ}) that a zonotope is a parallelotope if and only if summed segments are parallel to vectors of a unimodular set. McMullen was the first who discovered this fact in \cite{McM}. But he did not use the notion of unimodularity. Recall, that a set of vectors $X$ is {\em unimodular} if for any basic subset $B\subseteq X$ all vectors of $X$ have integer coordinates in $B$. An important unimodular set is the root system ${\mathbb A}_n$ for each dimension $n$.

Voronoi defines an {\em L-type} of a parallelotope which, in modern terms, is an isomorphism class of its face lattice. For brevity sake, we call L-type simply {\em type}.
Besides, Voronoi defined a type domain of a parallelotope as a set of parallelotopes of the same type. In particular, Voronoi distinguished a {\em principal} type domain. It is known and we show this in our paper, that parallelotopes of principal type domain are zonotopes related to the unimodular system ${\mathbb A}_n$.

We show that parallelotopes of closure of principal domain can be, up to an affine equivalence, described by the following system of inequalities
\begin{equation}
\label{PbS}
P(\beta)=\left\{x\in{\mathbb R}^N:\sum_{i\in S}x_i\le\beta(S) \mbox{  for all  }S \subset N,\mbox{  }\sum_{i\in N}x_i=0\right\},
\end{equation}
where $N=\{0,1,2,...,n\}$  is a finite set, and $\beta(S)$ is a submodular cut set function. A definition of submodular functions is given in Section 4.

The main property of the system (\ref{PbS}) is that it describes all parallelotopes of the closed principal domain. In other words each parallelotope of this domain can be obtained from any other by parallel shifts of faces defining hyperplanes.

Note  that polytopes of the form (\ref{PbS}) are called  {\em base polyhedra}. These polyhedra have integer vertices if values $\beta(S)$ are integer for all $S\subseteq N$. Submodular polytopes of type $P(\beta)$ are studied in many papers, see, for example, \cite{DK,Ed,Fu,Gr1} and others. They have numerous interesting applications, one of which we consider in this paper.

\section{A unimodular set ${\mathbb A}_n$}

Let ${\mathcal E}=\{e_i:i\in N\}$ be an orthonormal basis of the space ${\mathbb R}^{n+1}$. (Recall that $N=\{0,1,2,...,n\}$). For $S\subseteq N$, set $e_S=\sum_{i\in S}e_i$, and let ${\mathcal Q}_N=\{\pm e_S: S\subseteq N, S\not=\emptyset\}$.

The set ${\mathbb A}_n=\{\pm(e_i-e_j), i,j\in N, i<j\}$ is the well known root system lying in the hyperplane $\langle e_N,x\rangle=\sum_{i\in N}x_i=0$ of ${\mathbb R}^N$, where $\langle a,b\rangle$ denotes scalar product of vectors $a,b\in{\mathbb R}^N$. The assertion of the following Lemma~\ref{lt} is obvious.
\begin{lem}
\label{lt}
Let $t\in{\mathbb A}_n$, $t=e_i-e_j$. Then $\langle e_S,t\rangle=0$, if and only if either $\{ij\}\subseteq S$ or $\{ij\}\cap S=\emptyset$. Otherwise $\langle e_S,t\rangle\in\{\pm 1\}$.
\end{lem}

It is well known that the set ${\mathbb A}_n$ is unimodular, i.e. in any basic subset of ${\mathbb A}_n$ all its vectors have integer coordinates. Moreover, ${\mathbb A}_n$ is {\em graphic} unimodular system related to the complete graph $K_{n+1}$. See, for example, \cite{DK}.

Let $N$ be the set of vertices of $K_{n+1}$. Denote an edge of $K_{n+1}$ with vertices $i$ and $j$ by $(i,j)$.

There is a one-to-one correspondence between edges of $K_{n+1}$ and pairs of opposite vectors of ${\mathbb A}_n$. Namely, $(i,j)\Leftrightarrow \pm(e_i-e_j)$ for all $i,j\in V$, $i\not=j$. For an edge $u=(i,j)$ denote the vector $e_i-e_j$ by $t(u)=t(i,j)$, and, for $t\in{\mathbb A}_n$, denote the corresponding edge by $u(t)$.

\section{Perfect=principal domain ${\mathcal P}(A_n)$}

Roots of the set ${\mathbb A}_n$ generate integrally the root lattice $A_n$. By this definition, the lattice $A_n$ lies in the hyperplane $\langle e_N,x\rangle=0$. Points of $A_n$ are given in the orthonormal basis $\mathcal E$. The roots $t\in{\mathbb A}_n$ are {\em minimal} vectors of the root lattice $A_n$. Since there are, up to signs, $\frac{1}{2}n(n+1)$ minimal vectors, the lattice $A_n$ is {\em perfect}. Each root $t=e_i-e_j$ gives a quadratic form $\langle t,x\rangle^2= (x_i-x_j)^2$ of rank 1. These $\frac{1}{2}n(n+1)$ quadratic forms span extreme rays of the perfect domain ${\mathcal P}(A_n)\subset{\mathbb R}^{\frac{1}{2}n(n+1)}$ of quadratic forms. It coincides with the {\em principal} type domain, defined by Voronoi. Each quadratic form $f\in{\mathcal P}(A_n)$ has the following representation
\[f(x)=\sum_{t\in{\mathbb A}_n}b_t\langle t,x\rangle^2=\sum_{0\le i<j\le n}b_{ij}(x_i-x_j)^2, \]
where $x\in{\mathbb R}^N$, $\sum_{i\in N}x_i=0$ and $b_t=b_{ij}=b_{ji}$ is a coefficient related to the vector $t=e_i-e_j$.

Call a subset $X\subseteq{\mathbb A}_n$ {\em asymmetric} if it contains at most one vector from each pair $\{\pm t\}$ of opposite vectors. For each asymmetric subset $X\subset{\mathbb A}_n$, the forms $\langle t,x\rangle^2$ for all $t\in X$ generate a face of the principal domain ${\mathcal P}(A_n)$. Hence quadratic forms lying in the closure of ${\mathcal P}(A_n)$ are as follows
\begin{equation}
\label{fx}
f(x)=\sum_{t\in X}b_t\langle t,x\rangle^2.
\end{equation}
The Gram matrix of the form $f$ is
\begin{equation}
\label{DbX}
D_b(X)=\sum_{t\in X}b_ttt^T,
\end{equation}
where $t^T$ is transpose of the column vector $t$. The matrix $D_b(X)$ is a weighted sum of rank 1 matrices $tt^T$.

It is known, see, for example, \cite{RB,DV}, that the Voronoi polytope $P_V(f)$ of a quadratic form $f$ of an type domain is affinely equivalent to a Minkowski sum of Voronoi polytopes of forms lying on extreme rays of closure of this type domain. For each form $\langle t,x\rangle^2$ of rank 1 lying on extreme rays of ${\mathcal P}(A_n)$ its Voronoi polytope is the following segment
\begin{equation}
\label{zt}
z(t)=\left\{x\in{\mathbb R}^n:x=\lambda t, -1\le\lambda\le 1\right\}.
\end{equation}
This implies that, for $f$ defined in (\ref{fx}), $P_V(f)\simeq Z_b(X)$, where $\simeq$ denotes affine equivalence and
\[Z_b(X)=\sum_{t\in X}b_tz(t). \]
Obviously, $Z_b(X)$ is invariant with respect to changes vectors $t\in X$ by opposite vectors $-t$, and $z(t)$ is a one-dimensional zonotope.

\section{Graphs $G(X)$}
Parallelotopes of closure of the principal domain ${\mathcal P}(A_n)$ closely related to subgraphs $G\subseteq K_{n+1}$. In works of Russian authors, related graphs are called "symbols" (see, for example, \cite{Bo}, where parallelotopes of closure of principal domain are considered in detail).

In what follows, we consider only asymmetric sets of vectors.  Only  undirected graphs $G$ with sets of edges $U(G)$ are needed below for us. Hence, for $t\in{\mathbb A}_n$, $u(-t)=u(t)$.

Let $X\subset{\mathbb A}_n$ and let $U(X)=\{u:t(u)\in X\}$. The set of edges $U(X)$ generates a subgraph $G(X)\subseteq K_{n+1}$. We suppose that $N$ is vertex set of $G(X)$. If a vertex $v\in N$ is incident to no edge of $G(X)$, then it forms a one-vertex component $\{v\}$ of $G(X)$.

Similarly, let $G$ be a subgraph of $K_{n+1}$. Then $G$ determines a set $X(G)\subseteq{\mathbb A}_n$ such that $t\in X(G)$ if and only if $t=t(u)$ for some edge $u\in U(G)$.

Recall that rank of a subgraph $G\subseteq K_{n+1}$  is $rkG=n+1-c(G)$, where $c(G)$ is number of connected components of $G$. For $X\subseteq{\mathbb A}_n$, denote by dim$X$ dimension of a subspace of ${\mathbb R}^N$ spanned by vectors of $X$. Graph Theory gives the following equality: $rkG(X)={\rm dim}X$.

For $S\subseteq N$ denote by $U(S)\subseteq U(G)$ the set of edges of a subgraph $G(X)\subseteq K_{n+1}$ with one end in $S$ and another end in $\overline S$, where ${\overline S}=N-S$. Note that $U(S)=U({\overline S})$. Say that the set $S$ {\em cuts} the edges of the set $U(S)$. This set of edges is called a {\em cut} of the graph $G(X)$. Let $X_S=\{t\in X:u(t)\not\in U(S)\}$.  If $rkG(X_S)=rkG(X)-1$, then the set $U(S)$ is a {\em minimal} cut of the graph $G(X)$.

\section{Submodular cut set functions}
Let each edge $u\in U(K_{n+1})$ has a non-negative weight $b_u$. On the set $2^N$ of all subsets $S\subseteq N$ of the set of vertices of the graph $K_{n+1}$, define a non-negative set function  $\beta:2^N\to {\mathbb R}$ as follows:
\begin{equation}
\label{bS}
\beta(S)=\sum_{u\in U(S)}b_u.
\end{equation}
According to its definition, this function is called {\em cut function}. Since $U(S)$ is symmetric, $\beta(S)$ is also symmetric, i.e. $\beta({\overline S})=\beta(S)$. Hence one can consider only subsets $S\subseteq N$ such that $0\not\in S$. Besides, $\beta(\emptyset)=\beta(N)=0$.

It is well known, and this can be easily verified, that the function $\beta(S)$ is {\em submodular}, i.e. it satisfies the following inequalities
\[\beta(S)+\beta(T)\ge\beta(S\cap T)+\beta(S\cup T)\mbox{   for all  } S,T\subseteq N. \]
If these inequalities hold as equalities for all $S,T\subseteq N$, then the set function $\beta$ is called {\em modular}. For example, vectors $e_S$ considered as values of a vector set function form a modular vector set function.

The set of all cut functions forms a subcone of the cone of all  submodular functions taking zero value on empty set. The function $\beta$ is a point of the space ${\mathbb R}^{2^n}$. By definition, each cut function is uniquely determined by $\frac{1}{2}n(n+1)$ parameters $b_u$ for all $\frac{1}{2}n(n+1)$ edges $u$ of the graph $K_{n+1}$. Hence the cut cone lies in a space of dimension $\frac{1}{2}n(n+1)$.

Let
\[\delta_u(S)=\left\{
\begin{array}{ll}
1 &\mbox{if  }u\in U(S), \\
0 &\mbox{otherwise}.
\end{array}\right. \]

Recall that we relate to each vector $t\in{\mathbb A}_n$ an edge $u(t)\in U(K_{n+1})$ such that $u(t)=u(-t)$. Lemma~\ref{lt} can be stated as
\begin{lem}
\label{ds}
For any $t\in{\mathbb A}_n$ and $S\subseteq N$, the following equality holds
\[\delta_{u(t)}(S)=\langle e_S,t\rangle^2. \]
\end{lem}

Using the function $\delta_u(S)$ and Lemma~\ref{ds}, we can rewrite the equalities (\ref{bS}) for a graph $G\subseteq K_{n+1}$ in the form
\begin{equation}
\label{uS}
\beta_G(S)=\sum_{u\in U(G)}b_u\delta_u(S)=\sum_{t\in X}b_{u(t)} \langle e_S,t\rangle^2.
\end{equation}

The first equality in (\ref{uS}) shows that the functions $\delta_u$ are extreme rays of the cone of cut submodular functions. This gives an extreme ray description of this cone. But it is useful to know a facet description of the cone.

It is shown in \cite{Gr} that, for $n\ge 4$, a value $\beta(S)$ of a cut function on any subset $S\subseteq N$, where $|N|=n+1$, can be expressed through its values $\beta(\{ij\})=\beta(ij)$ on two-elements subsets $\{ij\}$. Here and below, for simplicity sake, for a set $S=\{ij...k\}$, we set $\beta(\{ij...k\})=\beta(ij...k)$.

Facets of the cut cone are defined by the following inequalities
\begin{equation}
\label{face}
\beta(i)+\beta(j)-\beta(ij)\ge  0 \mbox{   for all  }i,j\in N, i\not=j.
\end{equation}
A description of the linear space, where the cone of cut functions lies is given in \cite{Gr} and together with these inequalities it gives a facet description of the cone of cut set functions.

For $G\subseteq K_{n+1}$, denote by ${\mathcal C}(G)$ the cone of all cut functions defined on the graph $G$. Functions of this cone are defined in (\ref{uS}), where $b_u>0$ only for $u\in U(G)$. Hence the cone ${\mathcal C}(G)$ is open. Closure of this cone is defined by the conditions  $b_u=0$ for some $u\in U(G)$. For each $u\not\in U(G)$, $u=(i,j)$, the equality $\beta(ij)=\beta(i)+\beta(j)$ holds. Extreme rays of closure ${\mathcal C}(G)$ spanned by functions $\delta_u$ for $u\in U(G)$. This cone is simplicial, and has dimension $|U(G)|$. It is a $|U(G)|$-dimensional face of the whole cone of cut submodular functions ${\mathcal C}(K_{n+1})$.

Let $u=(i,j)$. Then definition of the function $\beta$ implies the following equality
\[\beta(i)+\beta(j)-\beta(ij)=2b_u. \]
This means that the set function $\beta$ uniquely determines the parameters $b_u$ for all $u\in U(K_{n+1})$. In particular, it determines uniquely a graph edges of which support parameters $b_u$.

\section{Submodular base polytopes}

Consider a polytope $P(\beta)$ related to a submodular function $\beta$:
\begin{equation}
\label{PSa}
P(\beta)=\left\{x\in{\mathbb R}^N:\langle e_S,x\rangle\le \beta(S) \mbox{  for all  }S\subset N, \mbox{  and  }\langle e_N,x\rangle=0\right\}.
\end{equation}
It is called {\em base polytope}. This polytope was intensely studied for its link with several integer optimization problems.

It is worth to note that if $\beta({\overline S})=\beta(S)$, then the system (\ref{PSa}) is equivalent to the following symmetric system of inequalities
\begin{equation}
\label{sym}
P(\beta)=\left\{x\in{\mathbb R}^n:-\beta(S)\le\langle e_S,x\rangle\le \beta(S) \mbox{  for all  }S\subseteq N-\{0\}\right\}.
\end{equation}
In fact, subtracting the equality $\langle e_N,x\rangle=0$ from the inequality $\langle e_{\overline S},x\rangle\le\beta({\overline S})=\beta(S)$ for each $S\subseteq N-\{0\}$, we obtain the system \ref{sym}. Besides this shows that the inequalities for $S$ and $\overline S$ define opposite faces.

For $\alpha\in{\mathbb R}$ and $p\in{\mathbb R}^N$, define a hyperplane
\begin{equation}
\label{hyp}
H(\alpha,p)=\{x\in{\mathbb R}^N:\langle p,x\rangle=\alpha\}.
\end{equation}

{\bf Remark}. One needs take in attention that the polytope $P(\beta)$ lies in the $n$-dimensional hyperplane $H(0,e_N)$. But facet vectors $e_S$ are defined as vectors of the space ${\mathbb R}^{n+1}$. Hence "real" facet vectors of $P(\beta)$ are projections of $e_S$ onto the hyperplane $H(0,e_N)$. In the description (\ref{sym}), one can consider coordinates $x_i$ in any basis of ${\mathbb R}^n$, in particular, in an orthonormal basis.

\vspace{2mm}
Base polytopes have many nice properties. The main property is that the set of base polytopes is closed under Minkowski sums. Obviously, if $\beta_1,\beta_2$ are submodular functions, then so is their sum $\beta_1+\beta_2$. Hence the following implication is true
\begin{equation}
\label{sum}
\beta=\beta_1+\beta_2 \Rightarrow  P(\beta)=P(\beta_1)+P(\beta_2).
\end{equation}

Another property of the polytope $P(\beta)$ is as follows (see, for example, \cite{Fu,Gr1}).
For $x\in P(\beta)$, let
\[{\mathcal S}(x)=\{S\subseteq N:\langle e_S,x\rangle=\beta(S)\}.\]
\begin{lem}
\label{lat}
For any submodular function $\beta$ and any $x\in P(\beta)$, the set ${\mathcal S}(x)$ is a sublattice of the Boolean $2^N$, i.e.
\[S,T\in{\mathcal S}(x) \mbox{  implies  }S\cap T,S\cup T\in{\mathcal S}(x).\]
\end{lem}
{\bf Proof}. If $S,T\in{\mathcal S}(x)$, then $\langle e_S,x\rangle=\beta(S)$, $\langle e_T,x\rangle=\beta(T)$. Recall that $e_S$ is a modular vector set function. Hence
\[\beta(S)+\beta(T)=\langle e_S,x\rangle+\langle e_T,x\rangle=
\langle e_{S\cap T},x\rangle+\langle e_{S\cup T},x\rangle\le\beta(S\cap T)+\beta(S\cup T). \]
Submodularity of the function $\beta$ implies that this is equality. This is possible only if $S\cap T,S\cup T\in{\mathcal S}(x)$. \hfill $\Box$

\vspace{2mm}
The following important theorem was proved firstly by J. Edmonds in \cite{Ed}, see also \cite{Fu}, \cite{Gr1}.
\begin{thm}
\label{ext}
A point $x$ is a vertex of $P(\beta)$ if and only if
\begin{equation}
\label{chain}
{\mathcal S}(x)=\{S_i:0\le i\le n, \mbox{  with  }|S_i|=i \mbox{ and  } S_i\subset S_{i+1}\mbox{ for }i<n\}
\end{equation}
is a chain of nested sets, where we set $S_0=\emptyset, S_{n+1}=N$. In this case, coordinates $x_i$ of the vertex $x$ are as follows
\begin{equation}
\label{xi}
x_i=\beta(S_i)-\beta(S_{i-1}), \mbox{  }1\le i\le n+1.
\end{equation}
\end{thm}

Theorem~\ref{ext} asserts that to each vertex of the polytope $P(\beta)$ there corresponds a full order $o$ on the set $N$ (or a permutation of the set $N$). Write $i<^oj$ if $i\in N$ stays in this order early than $j\in N$.

\section{Zonotopes}
Consider the base polytope $P(\beta)$ for a cut submodular function $\beta$. Coordinates of a vertex $x^o\in P(\beta)$ related to an order $o$ for the cut function $\beta$ are as follows.
\begin{equation}
\label{bij}
x_i^o=\sum_{j\in N:i<^oj}b_{ij}-\sum_{j\in N:j<^oi}b_{ij},
\end{equation}
where $b_{ij}=b_u$ if $u=(i,j)$. Note that several orders may correspond to the same vertex if $b_{ij}=0$ for some edges $(i,j)$.

Consider a graph $G\subseteq K_{n+1}$, and related to it cut submodular function $\beta_G$ defined in (\ref{uS}). Since $\beta_G=\sum_{u\in U(G)}b_u\delta_u$, by (\ref{sum}), we have
\[P(\beta_G)=\sum_{u\in U(G)}b_uP(\delta_u). \]
At first, we find vertices of $P(\delta_u)$ for $u=(i,j)$. So $b_u=0$ if $u\not=(i,j)$ and $b_{ij}=1$. By (\ref{bij}), all orders $o$, where $i<^oj$, correspond to one vertex with coordinates
\[x_i=1, \mbox{  }x_j=-1, \mbox{  }x_k=0\mbox{  for  }k\in N-\{i,j\}.\]
Similarly, all orders $o$, where $i>^oj$, correspond to a vertex with coordinates
\[x_i=-1, \mbox{  }x_j=1, \mbox{  }x_k=0\mbox{  for  }k\in N-\{i,j\}.\]
Therefore, the polytope $P(\delta_{(i,j)})$ is the segment $z(t)$ that is defined in (\ref{zt}). It is symmetric with respect to origin and is parallel to the vector $t=e_i-e_j\in{\mathbb A}_n$.

Recall that each $t$ corresponds to an edge $u=u(t)\in U(K_{n+1})$ with weight $b_u$. Set $b_t\equiv b_{u(t)}$.

So, we obtain
\[P(\delta_u)=z(t(u)),\mbox{  and  }P(\beta_G)= \sum_{t\in X} b_tz(t)=Z_b(X), \]
where $X=X(G)\subset{\mathbb A}_n$ is a set of roots related to edges $u\in U(G)$.

It is well known, see, for example, \cite{BZ}, and below, that $Z_b(K_{n+1})$ is a {\em permutohedron} if $b_u=b_0$ for all $u\in U(K_{n+1})$.

If $G\not=K_{n+1}$, then some of inequalities of the system (\ref{PSa}) are superfluous. It is sufficient to consider inequalities related to sets $S$ such that $U(S)$ is a minimal cut, when $e_S$ are facet vectors. There are different types of families of sets such that $e_S$ are facet vectors of $Z_b(G)$. For example, {\em intersecting}, {\em crossing}, {\em laminar}, {\em nested} and other families. Each of these families determines a submodular polytope $P(\beta)$ with some special properties, see \cite{Fu}.

The expression (\ref{PSa}) shows that any parallelotope of closure of principal domain may be obtained from any other by shifting face defining hyperplanes $H(\beta(S),e_S)$. We show below that the hyperplane $H(\beta(S),e_S)$ supports $P(\beta)$ at a face for each $S\subseteq N$.

\begin{prop}
\label{sup}
The hyperplane $H(\alpha_S,e_S)$ supports $P(\beta_G)=Z_b(X)$ at a face, for all $S\subseteq N$.
\end{prop}
{\bf Proof}. Obviously, for each $S\subseteq N$, there is some $\alpha\ge 0$ such that  the hyperplane $H(\alpha,e_S)$ supports a face $F_S$ of the zonotope $Z_b(X)$. Let $X_S=\{t\in X:\langle e_S,t\rangle=0\}$. By Lemma~\ref{lt}, $\langle e_S,t\rangle\in\{\pm 1\}$ for all $t\in X-X_S$. Hence (see, for example, \cite{DG1,McM1}), the center of the face $F_S$ is an end-point of the vector
\[q_S=\sum_{t\in X-X_S}b_t\langle  e_S,t\rangle t. \]
Since $\langle e_S,t\rangle=0$ for $t\in X_S$, we can sum over all $t\in X$. Obviously, for any point $x\in F_S$, and, in particular, for the center $q_S$ of the face $F_S$, the equality $\langle e_S,x\rangle =\alpha$ holds. Hence, using (\ref{uS}), we have
\[\alpha=\langle e_S,q_S\rangle=\sum_{t\in X}b_t\langle e_S,t\rangle^2=\beta_G(S). \]
\hfill $\Box$

\vspace{2mm}
Recall that each $(n-2)$-face $F$ of a parallelotope generates a 4- or 6-belt of facets containing $(n-2)$-faces that are parallel to $F$. Three, up to signs, facet vectors of a 6-belt defined by the face $F$ lie in a 2-plane, that is orthogonal to the $(n-2)$-face $F$. Hence, these three vectors are linearly dependent. For $Z_b(G)$, dependencies between three facet vectors are necessarily as follows
\begin{equation}
\label{6b}
e_S+e_T=e_{S\cup T},\mbox{  where  }S\cap T=\emptyset.
\end{equation}
Besides the three cuts $U(S), U(T), U(S\cup T)$ should be minimal. A particular but important case of equalities (\ref{6b}) when $S=\{i\}$ and $T=\{j\}$ are one-element sets, and $S\cup T=\{(ij)\}$ is an edge of the graph $G$.

\section{A type domain of the parallelotope $Z_b(G)$}

Let $\beta_0$ be a cut submodular function. Consider the polytope $P(\beta_0)$ defined in (\ref{PSa}). A set of all cut functions $\beta$ such that the polytope $P(\beta)$ has type of $P(\beta_0)$ consists of disjoint open connected components. Define {\em type domain} of the polytope $P(\beta_0)$ as a connected set containing $\beta_0$ and all cut submodular set functions $\beta$ determining polytopes $P(\beta)$ of type of the polytope $P(\beta_0)$.

\begin{thm}
\label{typ}
For a graph $G\subseteq K_{n+1}$, type domain of the parallelotope $Z_b(G)$ is the open cone ${\mathcal C}(G)$ of cut submodular set functions defined on the graph $G$. It is in one-to-one correspondence with a face of closure of principal domain of quadratic forms.
\end{thm}
{\bf Proof}. The zonotope $Z_b(G)$ depends of the set of parameters  $\{b_u>0: u\in U(G)\}$. It is obvious that changes of $b_u$ in this set correspond to homotheties in directions $t(u)$. These homotheties do not change the type of $Z_b(G)$. Obviously, these homotheties correspond to changes of cut submodular functions in the cone ${\mathcal C}(G)$.

Recall that the matrix $D_b(X)$ from (\ref{DbX}) generates a quadratic form $f=x^TD_b(X)x$ whose Voronoi polytope $P_V(f)$ is affinely equivalent to $Z_{\beta}(X)$. Consider the cut submodular function $\beta_X$ related to a set $X\subset{\mathbb A}_n$. Using (\ref{uS}), we have
\[\beta_X(S)=e_S^T(\sum_{t\in X}b_t(tt^T))e_S=e_S^TD_b(X)e_S.\]
We see that there is a one-to-one correspondence between cut set functions $\beta_X$ and quadratic forms determined by matrices $D_b(X)$, i.e. the domain of the Voronoi polytope $P_V(f)$, can be identified with the cone of cut submodular set functions. \hfill $\Box$

\vspace{2mm}
In \cite{BG} a {\em non-rigidity degree} of a lattice was defined as dimension of the type domain containing a form of this lattice. Hence it is natural to define a non-rigidity degree of a parallelotope $P$ and denote it nrd$P$ as dimension of its type domain. In other words, nrd$P$ is a number of independent parameters determining the parallelotope $P$.

Let a graph $G$ has $m=|U(G)|$ edges. Since the cone ${\mathcal C}(G)$ is simplicial, and, for $u\in U(G)$, the functions $\delta_u$ are its extreme rays, dimension of ${\mathcal C}(G)$ equals $m$. Hence nrd$Z_{\beta}(G)=m$. This implies that the system (\ref{PSa}) should have $m$ independent parameters. Since $\beta({\overline S})=\beta(S)$ and $\beta(N)=0$, this system of inequalities has $2^n-1$ parameters $\beta(S)$. By \cite{Gr}, parameters $\beta(S)$, for $|S|\not=2$, are represented through $\frac{1}{2}n(n+1)$ parameters $\beta(ij)$. If $G\not=K_{n+1}$, then $m<\frac{1}{2}n(n+1)$. Hence additional equalities should be between parameters $\beta(ij)$. Such dependencies provide sets $S$ that induce non-minimal cuts.

In fact, if $U(S)$ is non-minimal, then there are at least two sets $S_1$ and $S_2$ such that $S=S_1\cup S_2$, $S_1\cap S_2=\emptyset$ and $U(S)=U(S_1)+U(S_2)$, where sum denote disjoint union. In this case, the expression (\ref{bS}) gives  $\beta(S)=\beta(S_1)+\beta(S_2)$. We call such equalities {\em simple}.

Suppose that $p\in{\mathcal Q}_N$ is a facet vector of $Z_b(X)$. Then, for $D_b(X)$ from (\ref{DbX}), the vector
\[q_p=\sum_{t\in X}b_tt\langle t,p\rangle=D_b(X)p \]
is a minimal lattice vector, and we have linear map $q_p=D_{\beta}(X)p$ for all facet vectors.  This proves Voronoi conjecture for zonotopes (see, for the first proof \cite{Er}, and \cite{DG1} for this proof).

\section{Primitive polytopes $P(\beta)$}

An $n$-polytope  is called {\em primitive} if each its $k$-face (i.e. a face of dimension $k$) is contained exactly in $n-k$ facets.

\begin{lem}
\label{all}
If $P(\beta)$ is a primitive parallelotope, then vectors $e_S$ are facet vectors for all $S\subset N$, $S\not=\emptyset, N$.
\end{lem}
{\bf Proof}. A primitive $n$-parallelotope has $2(2^n-1)=2^{n+1}-2$ facets. All facets of $P(\beta)$ are determined by some of $2^{n+1}-2$ inequalities of the system (\ref{PSa}). If $P(\beta)$ is primitive, then all inequalities determine facets. \hfill $\Box$

\vspace{2mm}
Denote by $F_S$,  the facet of the primitive polytope $P(\beta)$ with facet vector $e_S$, $S\not=\emptyset,N$.

\begin{lem}
\label{all}
If $P(\beta)=Z_b(G)$ is a primitive parallelotope, then $G=K_{n+1}$.
\end{lem}
{\bf Proof}. If $G\not=K_{n+1}$, then for all $u\in U(K_{n+1})-U(G)$ we have $b_u=0$. Let $u=(i,j)$. By definition of the function $\beta$, we have $\beta(ij)=\beta(i)+\beta(j)$. Since $P(\beta)$ is primitive,
$F_{\{i\}}$, $F_{\{j\}}$ and $F_{\{ij\}}$ are facets. Let $x\in F=F_{\{i\}} \cap F_{\{ij\}}$ be a point of the $(n-2)$-face $F$. Then we have
\[\langle e_{\{j\}},x\rangle=\langle e_{\{ij\}},x\rangle- \langle e_{\{i\}},x\rangle=\beta(ij)-\beta(i)=\beta(j).\]
Since $x$ is an arbitrary point of $F$, we have that the $(n-2)$-face $F$ is contained in 3 facets. This contradicts the primitivity of $P(\beta)$. Hence $P(\beta)=Z_b(K_{n+1})$. \hfill $\Box$

\vspace{2mm}
It is not difficult to verify that the zonotope $Z_b(K_{n+1})$ has $(n+1)!$ vertices each corresponding to an ordering of $N$.

For the function $\beta_0$, defined by (\ref{bS}), where $b_u=a$ for all edges $u$, we have $\beta_0(S)=a|U(S)|=a|S|(n+1-|S|)$. Since in (\ref{xi}) $|S_i|=i$, coordinates of a vertex of the polytope $P(\beta_0)= Z_{\beta_0}(K_{n+1})$ are as follows
\[x_i=(n-2(i-1))a, \mbox{  }1\le i\le n+1. \]
So, vertices of $Z_{\beta_0}(K_{n+1})$ are in one-to-one correspondence with $(n+1)!$ permutations of the set $\{n,n-2,n-4,...,4-n,2-n,-n\}$, and this polytope is {\em permutahedron}.  Cf. Example 2.2.5 and Exercise 2.10 in \cite{BZ}.

\section{A parallelotope related to the root lattice $A_n$}

Let $G=C_{n+1}$ be a Hamiltonian $(n+1)$-circuit in the graph $K_{n+1}$. Any subset $X\subseteq X(C_{n+1})$ of cardinality $|X|<n+1$ is linearly independent. Hence each $k$-face, $0\le k\le n-1$, of the zonotope $Z_{\beta}(C_{n+1})$ is a $k$-dimensional parallelepiped. Such poytopes are called {\em cubical}.

We show that the zonotope $Z_b(C_{n+1})$ has the type of the Voronoi polytope $P_V(A_n)$ of the root lattice $A_n$. The classical root lattice $A_n$ is integrally generated by vectors of the root system  ${\mathbb A}_n$.

Roots are facet vectors of the Voronoi polytope $P_V(A_n)$, and the root system is a unimodular set. Moreover it is a maximal by inclusion unimodular set. Besides it is, up to a linear transformation, unique maximal by cardinality $\frac{1}{2} n(n+1)$ unimodular set of vectors.

We show that the set of $\frac{1}{2}n(n+1)$ facet vectors of  the zonotope $Z_b(C_{n+1})$ is unimodular. Recall that facet vectors of a zonotope $Z_b(G)$ are in one-to-one correspondence with minimal cuts $U(S)$ of the graph $G$. A minimal cut partitions the graph into two  connected components. Hence minimal cuts of the circuit $C_{n+1}$ partition this circuit into two connected chains of vertices.

From two chains we consider a chain that does not contain the vertex 0. Note that some chains consist of only one vertex. A facet vector related to a chain with the set of vertices $S\subseteq N$ is $e_S$. It can be considered as incidence vector of this chain.

It is not difficult to verify that there are $\frac{1}{2} n(n+1)$ subchains in the chain $C_{n+1}-\{0\}$. It is well known that incident vectors of subchains of a chain of length $k$ form a unimodular set isomorphic to ${\mathbb A}_k$. This implies that the set of all facet vectors of $Z_b(C_{n+1})$ and of $P_V(A_n)$ are affinely equivalent.

Since $C_{n+1}$ has $n+1$ edges, the type domain ${\mathcal C}(C_{n+1})$ of $Z_b(C_{n+1})$ has dimension $n+1$. Hence nrd$Z_{\beta}(C_{n+1})=n+1$, that is consistent with nrd${\mathbb A}_n=n+1$, what was proven in \cite{BG}.

\section{A parallelepiped $Z_b(T)$}
Let $T$ be a spanning tree of the graph $K_{n+1}$. The set of vectors $X(T)$ is a maximal independent set. Hence the zonotope $Z_b(X(T)) = Z_b(T)$ is a parallelepiped.

Each edge $u\in U(T)$ determines a minimal cut $U(S_u)=\{u\}$, where $S_u\not\ni 0$ is a set of vertices of a subtree of $T$ arisen after deleting the edge $u$. By (\ref{bS}), $\beta(S_u)=b_u$. Since $T$ has $n$ edges, the polytope $Z_b(T)$ has $n$ pairs of opposite facets. Let ${\mathcal S}(T)=\{S_u:u\in U(T)\}$. For $S,S'\in{\mathcal S}(T)$, the intersection $S\cap S'$ is equal to one of the three sets $\emptyset$, $S$ or $S'$. Such a family of sets is called {\em laminar} family. If $T$ is a Hamiltonian chain with 0 as an end vertex of the chain, then ${\mathcal S}(T)$ is a {\em nested} family that is a special case of laminar families (see, for example, \cite{Fu,DK}). It is not difficult to verify that all $n$ vectors $e_S$ for $S\in{\mathcal S}(T)$ form an independent set. This confirms once more that $Z_b(T)$ is a parallelepiped.

In this case, the system of inequalities in (\ref{sym}) is very special. It consists of $n$ pairs of inequalities. Each vertex of $Z_b(T)$ is a solution of $n$ equalities obtained from each pair of inequalities by substitution one of inequalities by equality. So, $Z_b(T)$ is a parallelepiped.

Let $T=K_{n,1}$ be a star, where all $n$ edges  are incident to a vertex, say, vertex 0. Then restriction of each function $\beta\in{\mathcal C}(K_{n,1})$ on sets $S\subseteq N-\{0\}$ is modular. Let $\beta_m$ one of these functions. Similarly as each modular set function taking zero value on empty set, the function $\beta_m$ is determined uniquely by its values $\beta_m(i)$ on one-element sets $\{i\}\subset N-\{0\}$ as follows $\beta_m(S)=\sum_{i\in S}\beta_m(i)$.

The polytope $P(\beta_m)$ in description (\ref{sym}) takes the form
\[P(\beta_m)=\{x\in{\mathbb R}^n:-\beta_m(i)\le x_i\le\beta_m(i)\mbox{   for all  }i\in N-\{0\}.\}\]
If we consider coordinates $x_i$ as coordinates of a point in an orthonormal basis of ${\mathbb R}^n$, then the polytope $P(\beta_m)$ is a rectangle parallelepiped whose edges have lengths $2\beta_m(i)= 2b_{(0,i)}$. It is an affine image of the zonotope $Z_b(K_{n,1})$.

\end{document}